\documentclass[12pt]{article}
\usepackage{amsfonts,amsmath,amsxtra}
\usepackage{amssymb}
\def\hybrid{\topmargin 0pt      \oddsidemargin 0pt
        \headheight 0pt \headsep 0pt
        \textwidth 16.5cm
        \textheight 23cm
        \marginparwidth 0.0in
        \parskip 5pt plus 1pt   \jot = 1.5ex}
\catcode`\@=11
\def\marginnote#1{}
\newcount\hour
\newcount\minute
\newtoks\amorpm
\hour=\time\divide\hour by60
\minute=\time{\multiply\hour by60 \global\advance\minute by-\hour}
\edef\standardtime{{\ifnum\hour<12 \global\amorpm={am}%
        \else\global\amorpm={pm}\advance\hour by-12 \fi
        \ifnum\hour=0 \hour=12 \fi
      \number\hour:\ifnum\minute<10 0\fi\number\minute\the\amorpm}}
\edef\militarytime{\number\hour:\ifnum\minute<10 0\fi\number\minute}

\def\draftlabel#1{{\@bsphack\if@filesw {\let\thepage\relax
   \xdef\@gtempa{\write\@auxout{\string
      \newlabel{#1}{{\@currentlabel}{\thepage}}}}}\@gtempa
   \if@nobreak \ifvmode\nobreak\fi\fi\fi\@esphack}
        \gdef\@eqnlabel{#1}}
\def\@eqnlabel{}
\def\@vacuum{}
\def\draftmarginnote#1{\marginpar{\raggedright\scriptsize\tt#1}}

\def\draft{\oddsidemargin -0.1truein
        \def\@oddfoot{\sl preliminary draft \hfil
        \rm\thepage\hfil\sl\today\quad\militarytime}
        \let\@evenfoot\@oddfoot \overfullrule 3pt
        \let\label=\draftlabel
        \let\marginnote=\draftmarginnote
\def\@eqnnum{{\rm (\theequation)}
\rlap{\kern\marginparsep\tt\@eqnlabel}%
\global\let\@eqnlabel\@vacuum}  }

\newcommand{\RR}{{\mathbb{R}}}
\newcommand{\CC}{{\mathbb{C}}}

\newcommand{\TT}{{\mathbb{T}}}

\newfont{\Bbbb}{msbm7 scaled 1\@ptsize00}
\newcommand{\zs}{\raise-1pt\hbox{$\mbox{\Bbbb Z}$}}

\font\teneufm=cmmib10 scaled 1\@ptsize00
\font\seveneufm=cmmib7 scaled 1\@ptsize00
\font\fiveeufm=cmmib5  
\def\bfit#1{{\textfont1=\teneufm\scriptfont1=\seveneufm
\scriptscriptfont1=\fiveeufm
\mathchoice{
\hbox{$\mathsurround=0pt\displaystyle#1$}}
{\mathsurround=0pt\hbox{$\textstyle#1$}}
{\hbox{$\mathsurround=0pt\scriptstyle#1$}}
{\hbox{$\mathsurround=0pt\scriptscriptstyle#1$}}}}

\font\sevenmsa=msam6 
\newfam\msafam
\textfont\msafam=\sevenmsa
\def\hexnumber@#1{\ifnum#1<10 \number#1\else
\ifnum#1=10 A\else\ifnum#1=11 B\else\ifnum#1=12 C\else
\ifnum#1=13 D\else\ifnum#1=14 E\else\ifnum#1=15 F\fi\fi\fi\fi\fi\fi\fi}
\def\msa@{\hexnumber@\msafam}
\def\llcorner{\delimiter"4\msa@78\msa@78 }
\def\lrcorner{\delimiter"5\msa@79\msa@79 }
\mathchardef\blacktriangleright="3\msa@49
\mathchardef\blacktriangleleft="3\msa@4A
\font\tenmsb=msbm10 scaled 1\@ptsize00
\newfam\msbfam
\textfont\msbfam=\tenmsb
\scriptfont\msbfam=\tenmsb
\def\msb@{\hexnumber@\msbfam}
\mathchardef\varkappa="0\msb@7B

\newdimen\linethick  \linethick=0.4pt
\newdimen\hboxitspace    \hboxitspace=5pt
\newdimen\vboxitspace    \vboxitspace=5pt

\def\fr#1{%
\be\new
\vcenter{
\hrule height\linethick
           \hbox{\vrule width\linethick
                 \kern\hboxitspace
                 \vbox{\kern\vboxitspace
                       \hbox{$\begin{array}{c}\displaystyle#1
          \end{array}$}%
                       \kern\vboxitspace}%
                 \kern\hboxitspace
                 \vrule width\linethick}%
           \hrule height\linethick}%
\ee}

\newdimen\Squaresize \Squaresize=14pt
\newdimen\Thickness \Thickness=0.5pt

\def\Square#1{\hbox{\vrule width \Thickness
   \vbox to \Squaresize{\hrule height \Thickness\vss
      \hbox to \Squaresize{\hss#1\hss}
   \vss\hrule height\Thickness}
\unskip\vrule width \Thickness}
\kern-\Thickness}

\def\Vsquare#1{\vbox{\Square{$#1$}}\kern-\Thickness}

\def\numberbysection{\@addtoreset{equation}{section}
        \def\theequation{\thesection.\arabic{equation}}}
\numberbysection

\renewcommand{\theequation}{\thesection.\arabic{equation}}
\def\titlepage{\@restonecolfalse\if@twocolumn\@restonecoltrue\onecolumn
     \else \newpage \fi \thispagestyle{empty}\c@page\z@
        \def\thefootnote{\fnsymbol{footnote}} }

\def\endtitlepage{\if@restonecol\twocolumn \else  \fi
        \def\thefootnote{\arabic{footnote}}
        \setcounter{footnote}{0}}  
\relax

\hybrid
\makeatletter
\newdimen\normalarrayskip            
\newdimen\minarrayskip               
\normalarrayskip\baselineskip
\minarrayskip\jot
\newif\ifold             \oldtrue            \def\new{\oldfalse}
\def\arraymode{\ifold\relax\else\displaystyle\fi}
\def\eqnumphantom{\phantom{(\theequation)}} 
\def\@arrayskip{\ifold\baselineskip\z@\lineskip\z@
     \else
     \baselineskip\minarrayskip\lineskip1\baselineskip\fi}

\def\@arrayclassz{\ifcase \@lastchclass \@acolampacol \or
\@ampacol \or \or \or \@addamp \or
   \@acolampacol \or \@firstampfalse \@acol \fi
\edef\@preamble{\@preamble
  \ifcase \@chnum
     \hfil$\relax\arraymode\@sharp$\hfil
     \or $\relax\arraymode\@sharp$\hfil
     \or \hfil$\relax\arraymode\@sharp$\fi}}

\def\@array[#1]#2{\setbox\@arstrutbox=\hbox{\vrule
     height\arraystretch \ht\strutbox
     depth\arraystretch \dp\strutbox
width\z@}\@mkpream{#2}\edef\@preamble{\halign \noexpand\@halignto
\bgroup \tabskip\z@ \@arstrut \@preamble \tabskip\z@ \cr}%
\let\@startpbox\@@startpbox \let\@endpbox\@@endpbox
  \if #1t\vtop \else \if#1b\vbox \else \vcenter \fi\fi
  \bgroup \let\par\relax
  \let\@sharp##\let\protect\relax
  \@arrayskip\@preamble}
\def\eqnarray{\stepcounter{equation}%
              \let\@currentlabel=\theequation
              \global\@eqnswtrue
              \global\@eqcnt\z@
              \tabskip\@centering              
              \let\\=\@eqncr
              $$%
            \halign to \displaywidth  \bgroup
             \eqnumphantom \@eqnsel
      \hskip\@centering                               
    $\displaystyle  \tabskip\z@ {##}$%
    &\global\@eqcnt\@ne \hskip 2\arraycolsep
         $ \displaystyle  \arraymode{##}$\hfil
    &\global\@eqcnt\tw@ \hskip 2\arraycolsep
         $\displaystyle\tabskip\z@{##}$\hfil
         \tabskip\@centering
    &{##}\tabskip\z@\cr}
\makeatother
\newtheorem{te}{Theorem}[section]
\newtheorem{de}{Definition}[section]
\newtheorem{prop}{Proposition}[section]           

\newcommand{\beq}[1]{\begin{equation}\label{#1}}
\newcommand\eeq{\end{equation}}
\newcommand\bqa{\begin{eqnarray}}
\newcommand\eqa{\end{eqnarray}}
\def\be{\begin{eqnarray}\new\begin{array}{cc}}
\def\ee{\end{array}\end{eqnarray}}

\def\beq{\begin{equation}}
\def\eeq{\end{equation}}
\def\bse{\begin{subequations}}                
\def\ese{\end{subequations}}
\def\bp{\begin{pmatrix}}
\def\ep{\end{pmatrix}}

\def\d{\partial}
\def\stack#1#2{\raise0.7pt\hbox{$\mathrel{\mathop{#2}\limits^{#1}}$}}
\def\tr{\triangleright}
\def\tl{\triangleleft}
\def\sem{\mathsurround=0pt \raise1pt
\hbox{$\scriptscriptstyle>\!\!$}\:\!\!\tl}
\def\mes{\mathsurround=0pt \tr\!\:\!\raise0.8pt
\hbox{$\scriptscriptstyle\!\!<$}\,}
\def\]{\mathsurround=0pt ]\raise-2pt\hbox{$_\ast$}}

\def\bgamma{{\bfit\gamma}}

\def\brho{{\bfit\rho}}

\def\la{\lambda}

\def\<{\langle}
\def\>{\rangle}
\def\ov{\overline}

\def\wh{\widehat}
\def\frak{\mathfrak}
\def\ca{{\cal A}}
\def\cb{{\cal B}}
\def\cc{{\cal C}}

\def\N{{\scriptscriptstyle N}}
\def\ts#1#2{{\textstyle\frac{#1}{#2}}}

\def\we{\raise-1pt\hbox{$\,\stackrel{\wedge}{,}\,$}}

\def\pr {\partial}

0
\begin{document}

\title{
\hfill{\normalsize ITEP-TH-1/04}\\ [10mm] \bf  Representation
theory and quantum integrability \footnote{Extended talk by the
third author at the satellite XIV ICMP workshop:
 "Infinite Dimensional Algebras and Quantum Integrable systems",
 July 21-25, 2003, Faro, Portugal }}

\author{A. Gerasimov\thanks{Institute for Theoretical and Experimental
Physics, Moscow, Russia}\,\,\thanks{Hamilton Mathematical
Institute at Trinity College,
Dublin, Ireland }%
\and S. Kharchev\thanks{Institute for Theoretical and Experimental
Physics, Moscow, Russia}%
\and D. Lebedev\thanks{Institute for Theoretical and Experimental
Physics, Moscow, Russia}}
\date{\phantom{.}}
\maketitle

\begin{abstract}
We describe  new constructions of the infinite-dimensional
representations of $U(\mathfrak{g})$ and $U_q(\mathfrak{g})$ for
$\mathfrak{g}$ being $\mathfrak{gl}(N)$ and $\mathfrak{sl}(N)$.
The application of these  constructions
to the quantum integrable theories of Toda type is discussed. With
the help of these infinite-dimensional representations we manage
to establish direct connection between group theoretical  approach
to the quantum integrability and Quantum Inverse Scattering Method
based on the representation theory of Yangian and its
generalizations.

In the case of $U_q(\mathfrak{g})$ the considered representation
is naturally supplied with the structure of
$U_q(\mathfrak{g})\otimes U_{\tilde
q}(\check{\mathfrak{g}})$-bimodule where $\check {\mathfrak{g}}$
is Langlands dual to $\mathfrak{g}$ and $\log q/2\pi i=-
(\log{\tilde q}/2\pi i)^{-1}$. This bimodule structure is a
manifestation of the Morita equivalence of  the
algebra and its dual.
\end{abstract}
\clearpage \newpage


\tableofcontents
\normalsize
\section{Introduction}

Since the early days of quantum mechanics,  the representation theory
plays the important role as a succinct tool to describe and
explicitly solve quantum theories. On the other hand, the problems
in quantum theories serve as a source of the new ideas in
representation theory. The most notable recent example is the
emergence of the notion of the quantum groups \cite{Dr0}, \cite{J} from
the algebraic formulation of the Quantum Inverse Scattering Method (QISM)
\cite{F}, \cite{KS}.

In this paper we describe several constructions in representation
theory of classical and quantum groups inspired by our studies of
the simple quantum  integrable models. Our starting point was the
desire to understand QISM in more standard representation theory
terms. It was a surprise that one should instead  invent the new
constructions in the representation theory. It is well known to
experts in integrable systems that  there exist distinguished
coordinates in which the description of the quantum system,
although non trivial, is drastically simplified. Multi-dimensional
nonlinear spectral problems are reduced to  the one dimensional
case and
the solution of the quantum theory may be constructed explicitly.
The recent  interest in this approach initiated by Sklyanin
\cite{Skl} leads to an explosion  of  activity connected with a
search of the explicit transforms to such coordinates. This
Separation of Variables method  (SoV) may be considered as an
extension/generalization of QISM.

It was shown in \cite{GKL}, \cite{GKL1} that the separation of variables
in the modern form has a clear group theoretic meaning and goes
back to the natural parameterization of the regular group elements
of the non-compact groups. This leads to the new construction of
the representations of the universal enveloping of the Lie algebras
 in terms of the difference  operators. A certain particular case of this
 construction leading to the finite dimensional representations
turns out to be  the well known Gelfand-Zetlin recursive
construction of the representations of  classical groups
\cite{GZ}, \cite{GG}, \cite{LP}.

The explicit connection with the QISM approach appears through the
closely related construction of  the representations of the
Yangian. Namely, let $\frak{g}$ be the Lie algebra $\frak{gl}(N)$
and let $Y(\frak{g})$ be the Yangian \cite{Dr2}. There is a natural
epimorphism ${\pi}_\N :Y(\frak{g})\rightarrow{\cal U}(\frak{g})$
compatible with the representation in terms of the difference
operators. The constructed representation of the Yangian turns out
to be a manifestation of the simultaneous existence of the
$RTT$-type realization and Drinfeld's "new" realization \cite{Dr1}. The
distinguished maximal commutative subalgebra of the Yangian is a
key ingredient in the explicit connection between these two
constructions and may be described by the set of the commuting
operators ${\cal A}_n(\lambda)\,,n=1,\ldots,N$. In the special
class of the representations ${\cal A}_n(\lambda)$ act as polynomials in
$\lambda$. The zeros of these polynomials provide the variables
appearing in SoV. On the other hand, these variables are exactly
the variables that appear in our generalization of Gelfand-Zetlin
construction.
 Note that  Drinfeld's "new" realization
is known for the Yangian of an arbitrary simple Lie algebra and
the generalized Gelfand-Zetlin
construction of the representation of $\frak{gl}(N)$ may be
naturally extended to the  case of the general simple Lie algebra.

There is a  further generalization of the construction providing
the explicit infinite-dimensional representation of the quantum
group $U_q(\mathfrak{g})$ for
$\mathfrak{g}=\mathfrak{gl}(N),\mathfrak{sl}(N)$.  However, a new
essential phenomenon  appears in this case. The
infinite-dimensional representations of $U_q(\mathfrak{g})$ are
naturally supplied with the structure of $U_q(\mathfrak{g})\otimes
U_{\tilde q}(\check{\mathfrak{g}})$-bimodule where $\check
{\mathfrak{g}}$ is Langlands dual to $\mathfrak{g}$ and
$\log q/2\pi i=-(\log{\tilde q}/2\pi i)^{-1}$ (for the
preliminary results in these direction see
\cite{Fad}, \cite{Fad1}, \cite{KLS}). We give explicit construction of
this bimodule structure and show the natural appearance of the
Langlands dual for $U_q(\mathfrak{g})$, when
$\mathfrak{g}=\mathfrak{gl}(N),\mathfrak{sl}(N)$. Note that
in this paper we essentially use  various rational forms of the
quantum groups (see \cite{Lu}, \cite{CP} and references therein). Thus,
the Langlands dual to the minimal ("adjoint") rational form of
$U_q(\mathfrak{sl}(N))$ turns out to be another form of
$U_{\tilde q}(\mathfrak{sl}(N))$ in accordance with the fact that the Langlands
dual to $PSL(N)$ is  $SL(N)$ for the  classical groups.
Further study shows that the appearance
of the Langlands dual in this construction is not accidental and
is a manifestation of the Morita equivalence
of certain two algebras
naturally associated with $U_q(\mathfrak{g})$ and $U_{\tilde
q}(\check{\mathfrak{g}})$  (compare with  \cite{Co},  \cite{Ri}).
This formulation  appears to be very
close to the initial problem of the classification of the irreducible
representations in the decomposition of the regular representation
of the group $G$ in terms of the dual group $\check{G}$ \cite{Vog}.
We are planning to discuss this construction in full details elsewhere.

The plan of the paper is as follows. In section 2 we review the
part of \cite{GKL} concerning  the analytic continuation of the
Gelfand-Zetlin construction to the case of Whittaker modules of
$U(\frak{gl}(N))$. An explicit description of a Whittaker module
is straightforward as soon as we have an explicit expression for the
cyclic Whittaker vector.
 In section 3 we generalize this construction to the  case  of quantum
 groups. The structure of $U_q(\mathfrak{g})\otimes
U_{\tilde q}(\check{\mathfrak{g}})$-bimodule and its  connection with
Langlands duality is shortly discussed.
In section 4 we demonstrate the connection with SoV and QISM and
describe the applications to the quantum integrable theories of
Toda type.

The results of these paper are based on the work  that was
done for the past few years by the ITEP group in Moscow and
 was reported in the series of papers \cite{KL0}, \cite{KL}, \cite{KL1},
\cite{KLS}, \cite{GKL}, \cite{GKL1}.

{\bf Acknowledgments:}
The research was partly supported by grants
CRDF RM1-2545-MO-03; INTAS 03-513350; grant 1999.2003.2
for support of scientific schools, and by grants
RFBR 03-02-17554 (A. Gerasimov, D. Lebedev), RFBR 03-02-17373 (S.
Kharchev).
We are grateful to M.Kont\-se\-vich, and M. Semenov-Tian-Shansky for
their  interest in this work and we are grateful to A. Rosly
for useful discussion.  D.L. is also  grateful
to IHES for warm hospitality and to  organizers and participants
of the ICMP workshop in Faro for the creation of the stimulating
atmosphere. D.L. thanks to 21 COE RIMS Research Project 2004:
Quantum Integrable Systems and Infinite Dimensional Algebras,
where the paper was finished, for support and warm hospitality.

\section{Whittaker modules in the Gelfand-Zetlin representation}

\subsection{The representation of $U(\frak{gl}(N))$}

 Let us  remind the construction of the paper \cite{GKL}, where
an analytical continuation
the Gelfand-Zetlin (GZ) theory  to
infinite-dimen\-sio\-nal representations of the universal enveloping
algebra ${ U}(\frak{gl}(N))$ was introduced.

 Let $\hat\TT_{\hbar}$ be an associative algebra generated
by $\hat{\gamma}_{nj}, \hat{\beta}_{nj}^{\pm 1}\;,n=1,\ldots,N-1;
1\leq j\leq n$, and
$\hat{\gamma}_{\N j}\;,1\leq j\leq N$, subject to the
relations
\be\label{ah}
[\hat{\gamma}_{nj}, \hat{\gamma}_{ml}]=
[\hat{\beta}_{nj},\hat{\beta}_{ml}]=0\,,\\

[\hat{\beta}_{ml},\hat{\gamma}_{nj}]=i\hbar\delta_{mn}\delta_{lj}
\hat{\beta}_{ml}\,.
\ee
\begin{te}
Let $E_{nm},\,n,m=1,\ldots,N$ be the generators of $\frak{gl}(N)$.
The following explicit expressions defines the embedding $\pi$:
 $\frak{gl}(N)\hookrightarrow \hat\TT_{\hbar}$
\bse\label{m1}
\be\label{m1a}
\hspace{2cm}
E_{nn}=\frac{1}{i\hbar}\Big(\sum_{j=1}^n\hat{\gamma}_{nj}-
\sum_{j=1}^{n-1}\hat{\gamma}_{n-1,j}\Big),\hspace{2cm}(n=1,\ldots,N),
\ee
\be\label{m1b}
E_{n,n+1}=-\frac{1}{i\hbar}\sum_{j=1}^n\frac{\prod\limits_{r=1}^{n+1}
(\hat{\gamma}_{nj}-\hat{\gamma}_{n+1,r}-\frac{i\hbar}{2})}
{\prod\limits_{s\neq j}(\hat{\gamma}_{nj}-\hat{\gamma}_{ns})}\,
\hat{\beta}^{-1}_{nj},
\hspace{0.5cm}
(n=1,\ldots,N-1),
\hspace{-2cm}
\ee
\be\label{m1c}
E_{n+1,n}=\frac{1}{i\hbar}\sum_{j=1}^n\frac{\prod\limits_{r=1}^{n-1}
(\hat{\gamma}_{nj}-\hat{\gamma}_{n-1,r}+\frac{i\hbar}{2})}
{\prod\limits_{s\neq j}(\hat{\gamma}_{nj}-\hat{\gamma}_{ns})}\,
\hat{\beta}_{nj},\hspace{0.5cm}(n=1,\ldots,N-1).
\hspace{-2cm}
\ee
\ese
\end{te}

Let us consider the following natural representation of the quantum torus
$\hat\TT_{\hbar}$:
\be\label{htori1}
\hspace{1cm}
\hat{\gamma}_{nj}=\gamma_{nj}\in\CC\,,\ \ \
n=1,\ldots,N\,,\ \  1\leq j\leq n\,,\\
\hspace{1.5cm}
\hat{\beta}_{nj}=e^{i\hbar\frac{\d}{\d\gamma_{nj}}}\,,\ \ \
n=1,\ldots,N-1\,,\ \ 1\leq j\leq n\;,
\ee
where $\gamma_{nj}$ and $e^{\pm i\hbar\frac{\d}{\d\gamma_{nj}}}$ are
considered as operators acting on the space $M$ of meromorphic functions
of complex variables $\gamma_{nj}, n=1,\ldots,N\!-\!1,1\leq j\leq n$. The
remaining
$\gamma_{\N1},\ldots,\gamma_{\N\N}$ are considered as constants. Thus the
complex vector
$\bgamma_\N=(\gamma_{\N1},\ldots,\gamma_{\N\N})$ plays the role of a label
which determines the above representation.

Let ${\cal Z}(\frak{gl}(n))$ be the center of ${ U}(\frak{gl}(n))$.
We say  that a  ${ U}(\frak{gl}(n))$-module $V$ admits an
infinitesimal character $\xi$ if there is a homomorphism
$\xi: {\cal Z}(\frak{gl}(n))\rightarrow \CC$ such that
$zv=\xi(z)v$ for all $z\!\in\!{\cal Z}(\frak{gl}(n)), v\in V.$
It is possible to show that the ${ U}(\frak{gl}(N))$-module $M$
defined above admits an infinitesimal character and
each  central element of ${ U}(\frak{gl}(n))$ acts on $M$  via
multiplication by a symmetric polynomial in the variables $\gamma_{nj}$
\cite{GKL}.

In the next subsection we calculate the explicit action of the central
elements on $M$ using the notion of Whittaker vectors.

\subsection{Whittaker modules}

Let us now  give a construction for  Whittaker modules
using the representation of ${ U}(\frak{gl}(N))$ described above.
We first recall some facts from \cite{Ko2}.

Let $n_+$ and $n_-$ be the subalgebras of $\frak{gl}(N)$
generated,  respectively,  by positive and negative root generators.
The homomorphisms (characters) $\!\chi_+\!\!:n_+\rightarrow\CC$,
$\chi_-\!\!:n_-\rightarrow\CC$
are uniquely determined by their values on the simple root generators,
and are called non-singular if
the (complex) numbers $\chi_+(E_{n,n+1})$ and
$\chi_-(E_{n+1,n})$ are non-zero for all $n=1,\ldots,N-1$.

Let $V$ be any ${ U}={ U}(\frak{gl}(N))$-module. Denote  the
action of $u\in{ U}$ on $v\in V$ by $uv$.
A vector $w\in V$ is called a Whittaker vector with respect to the
character $\chi_+$  if
\be\label{ewt}
\hspace{2cm}
E_{n,n+1}w=\chi_+(E_{n,n+1})w\,,
\hspace{1cm}(n=1,\ldots,N-1),
\ee
and an element $w'\in V'$ is called a Whittaker vector
with respect to the character $\chi_-$  if
\be\label{fwt}
\hspace{2cm} E_{n+1,n}w'=\chi_-(E_{n+1,n})w'\,,
\hspace{1cm} (n=1,\ldots,N-1).
\ee

A Whittaker vector is cyclic for $V$ if ${ U}w=V$,
and a ${ U}$-module is a Whittaker module
if it contains a cyclic Whittaker vector. The
${ U}$-modules $V$ and $V'$ are called dual
if there exists a non-degenerate pairing $\<.\,,.\>: V'\times V\to\CC$
such that
$\langle X v',v\rangle = -\langle v',Xv\rangle$ for all $v\in V,$
$v'\in V'$ and  $X\in\frak{gl}(N)$.

We proceed with explicit formulas for Whittaker vectors
corresponding to the representation given by (\ref{m1}), (\ref{htori1}).
\begin{prop}
The equations
\be\label{fww'}
E_{n+1,n}w_\N' =-i\hbar^{-1}w_\N' ,
\ee
\be\label{fww''}
E_{n,n+1}w_\N=-i\hbar^{-1}w_\N
\ee
for all $n=1,\ldots,N-1$, admit the solutions
\be\label{wv}
w_\N'=1 ,\\
w_\N=e^{-\frac{\pi}{\hbar}\sum\limits_{n=1}^{N-1}(n-1)
\sum\limits_{j=1}^n
\gamma_{nj}}\prod_{n=1}^{N-1} s_n(\bgamma_n,\bgamma_{n+1}),
\ee
where
\be\label{wv'}
s_{n}(\bgamma_n,\bgamma_{n+1})=\prod_{k=1}^n\prod_{m=1}^{n+1}
\hbar^{\frac{\gamma_{nk}-\gamma_{n+1,m}}{i\hbar}+\frac{1}{2}}\;
\Gamma\Big(\frac{\gamma_{nk}-\gamma_{n+1,m}}{i\hbar}+\frac{1}{2}\Big).
\ee
\end{prop}
 ( For the proof see \cite{GKL}.)

 The solutions (\ref{wv}), (\ref{wv'}) are not unique.
Indeed, the set of  Whittaker vectors is closed under the multiplication
by an arbitrary $i\hbar$-periodic function in the  variables $\gamma_{nj}$.
Hence, there are infinitely many invariant subspaces in $M$ corresponding
to infinitely many Whittaker vectors.

To construct irreducible submodules, let us introduce the Whittaker
modules $W$ and $W',$ generated cyclically by the Whittaker vectors
$w_\N$ and $w_\N',$ respectively. An explicit description of a Whittaker
module is straightforward as soon as we have an explicit expression for
the  Whittaker vectors.
Namely
let $\bfit m_n=(m_{n1},\ldots,m_{nn})$ be the set of non-negative integers.
The Whittaker module $W={ U}w_\N$ is spanned by the elements
\be\label{bas1}
w_{\bfit m_1,\ldots,\bfit m_{\N-1}}=\prod\limits_{n=1}^{N-1}
\prod\limits_{k=1}^{n}\sigma_{k}^{m_{nk}}(\bgamma_{n})w_\N,
\ee
where $\!\sigma_{k}(\bgamma_{n})\!$ is the elementary symmetric function of
the variables $\!\gamma_{n1},\ldots,\gamma_{nn}\!$ of order $\!k\!$:
\be\label{bas2}
\sigma_k(\bgamma_n)=\sum_{j_1<\ldots<j_k}
\gamma_{nj_1}\ldots\gamma_{nj_k}.
\ee
Similarly, the Whittaker module $W'={ U}w_\N'$ is spanned by the
polynomials
\be
w'_{\bfit m_1,\ldots,\bfit m_{\N-1}}=
 \prod\limits_{n=1}^{N-1}
\prod\limits_{k=1}^{n}\sigma_{k}^{m_{nk}}(\bgamma_{n}).
\ee
The Whittaker modules $W$ and $W'$ are irreducible.

Let us note that for any  subalgebra
${ U}(\frak{gl}(n))\subset{ U}(\frak{gl}(N))$, $2\leq n<N$,
the module over the ring of the polynomials in $\bgamma_{n}$ with the
basis $\prod\limits_{l=1}^{n-1} \prod\limits_{k=1}^{l}
\sigma_{k}^{m_{lk}}(\bgamma_{l})w_\N$ is a ${ U}(\frak{gl}(n))$
Whittaker module.
One can calculate the explicit form of the action of the central
elements of ${ U}(\frak{gl}(n))$ on the space $M$.
It is well known that the generating function $\ca_n(\la)$ of the
central elements  of ${ U}(\frak{gl}(n))$ (the Casimir
operators) can be represented as follows
\cite{HU}:
\be\label{cas1}
{\cal A}_n(\la)\\
\hspace{-8mm}=
\sum_{s\in S_n}{\rm sign}\,s
\Big[(\la-i\hbar\rho^{(n)}_1)\delta_{s(1),1}-
i\hbar E_{s(1),1}\Big]\ldots
\Big[(\la-i\hbar\rho^{(n)}_n)\delta_{s(n),n}-
i\hbar E_{s(n),n}\Big],\hspace{-0.8cm}
\ee
where $\rho^{(n)}_k=\frac{1}{2}(n-2k+1),\;k=1,\ldots, n$ and
the summation is over elements of the permutation group $S_n$.
It can be proved \cite{GKL} that
the operators (\ref{cas1}) have the following form  on a space
of meromorphic functions $M$:
\be\label{cas2}
\hspace{2cm}
\ca_n(\la)=\prod_{j=1}^n(\la-\gamma_{nj}),\hspace{1cm}(n=1,\ldots,N).
\ee

It  remains to construct a pairing between $W$ and $W'$, and to prove that
the Whittaker modules  $W$ and $W'$ are dual with respect to this pairing.
Let $\phi\in W'$ and $\psi\in W$. Define the pairing
$\langle.\,,.\rangle$: $W'\otimes W\rightarrow\CC$ by
\be\label{in1'}
\langle\phi,\psi\rangle=\int\limits_{\RR^{\frac{N(N-1)}{2}}}
\mu_0(\bgamma)\ov{\phi}(\bgamma)\,\psi(\bgamma)\,
\prod\limits_{\stackrel{\scriptstyle n=1}{j\leq n}}^{N-1}d\gamma_{nj}\,,
\ee
where
\be\label{in2'}
\mu_0(\bgamma)=\prod_{n=2}^{N-1}\prod\limits_{p<r}
(\gamma_{np}-\gamma_{nr})
(e^{\frac{2\pi\gamma_{np}}{\hbar}}-e^{\frac{2\pi\gamma_{nr}}{\hbar}}).
\ee
The integral (\ref{in1'}) converges absolutely (See for the proof \cite{GKL}).

 To construct the pair of the dual Whittaker modules,
we should restrict the label of
representation to the real values: $\bgamma_\N\in\RR^N$.
Then the Whittaker modules $W$ and $W'$ will be dual with
respect to the pairing defined by (\ref{in1'}). That is for any
$\phi\in W'$ and $\psi\in W,$ the generators $X\in\frak{gl}(N,\RR)$
possess the property
\be\label{inv5}
\langle\phi, X\psi\rangle\,=\,-\,\langle X\phi,\psi\rangle.
\ee
This property will be important in the derivation of the wave function
in Section 4.

\section{Construction of the representation of $U_{q}(\frak g)$}
In this section we outline an extension of our approach to the
quantum groups. We start with the construction of the embedding of
$U_{q}(\frak g)$ in the product of a commutative and
non-commutative tori. The skew field of fractions  constructed
from $U_{q}(\frak g)$ coincides with the functions on the product
of the tori invariant under the product of the symmetric groups.
Then we describe explicitly the structure of
$U_{q}(\frak g)\,\otimes U_{\tilde q}(\frak{\check{g}})$ bimodule.

\subsection{The  rational forms of  $U_{q}(\frak g)$}
We start with the definition  of the quantum groups following
\cite{Lu}, \cite{CP}.
Let $q$ be an indeterminate, and let $\CC (q)$ be the field of
rational functions of $q$ with coefficients in $\CC$.
The associative $\CC (q)$--algebra $U_{q}(\frak{gl}(N))$ is generated by the
elements $K^{\pm 1}_{nn}\,,E_{nm}\,;n\neq m\,;n,m=1,\ldots,N$
subjected to the relations\footnote{We will not use the Hopf structure of the
quantum groups and, therefore, we omit the comultiplication formulas in what
follows.
}:
\be\label{d1}

K_{nn}K_{nn}^{-1}=K_{nn}^{-1}K_{nn}=1\,,\ \ \
K_{nn}K_{mm}=K_{mm}K_{nn}\,,\\
K_{nn} E_{m,m+1}K^{-1}_{nn}=q^{\delta_{nm}-\delta_{n,m+1}}E_{m,m+1}\,,\\
K_{nn}E_{m+1,m}K_{nn}^{-1}=q^{\delta_{n,m+1}-\delta_{nm}}E_{m+1,m}\,,\\
E_{n,n+1}E_{m+1,m}- E_{m+1,m}E_{n,n+1}=
\delta_{nm}\frac{K_{nn}K_{n+1,n+1}^{-1}-K^{-1}_{nn}K_{n+1,n+1}}{q-q^{-1}}
\ee
together with quantum analogues of  Serre relations
\bse\label{d3}
\be\label{d3a}
\hspace{1cm}
E_{n,n+1}E_{m,m+1}-E_{m,m+1}E_{n,n+1}=0\,,\ \ \ \ m\neq n\pm1\,,\\
\hspace{-0.2cm}
E_{n,n+1}^2E_{n+1,n+2}-(q+q^{-1})E_{n,n+1}E_{n+1,n+2}E_{n,n+1}+
E_{n+1,n+2}E_{n,n+1}^2=0\,,\\
\hspace{-0.2cm}
E_{n+1,n+2}^2E_{n,n+1}-(q+q^{-1})E_{n+1,n+2}E_{n,n+1}E_{n+1,n+2}+
E_{n,n+1}E_{n+1,n+2}^2=0\,,
\ee
\be\label{d3b}
\hspace{1cm}
E_{n+1,n}E_{m+1,m}-E_{m+1,m}E_{n+1,n}=0\,,\ \ \ \ m\neq n\pm1\,,\\
\hspace{-0.2cm}
E_{n+1,n}^2E_{n+2,n+1}-(q+q^{-1})E_{n+1,n}E_{n+2,n+1}E_{n+1,n}+
E_{n+2,n+1}E_{n+1,n}^2=0\,,\\
\hspace{-0.2cm}
E_{n+2,n+1}^2E_{n+1,n}-(q+q^{-1})E_{n+2,n+1}E_{n+1,n}E_{n+2,n+1}+
E_{n+1,n}E_{n+2,n+1}^2=0\,.
\ee
\ese

We will also need the explicit description of the quantum group
$U_{q}(\frak{sl}(N))$. Note that there exist various rational
forms of quantum groups \cite{Lu}, \cite{CP}. In the case of
$U_{q}(\frak{sl}(N))$ we have the following description. Let
$a_{nm}=2\delta_{nm}-\delta_{n,m+1}-\delta_{n+1,n}$ be the Cartan matrix of
$\frak{sl}(N)$. Let $Q$ and $P$
denote the root lattice and the lattice of  weights, respectively.
The smallest rational form, the adjoint rational form
$U_{q}^{Q}(\frak{sl}(N))$, is the associative $\CC (q)$--algebra  with
generators $E_n , F_n$ and $ K_{n}^{\pm 1}\,, n=1,\ldots,N-1$, and
relations
\bse\label{sln}
\be\label{slna}
K_{n}K_{n}^{-1}=K_{n}^{-1}K_{n}=1\,,\ \ \ \ K_{n}K_{m}=K_{m}K_{n}\,,
\ee
\be\label{slnb}
K_{n}E_{m}K_{n}^{-1}=q^{a_{nm}}E_{m}\,,\ \ \ \
K_{n}F_{m}K_{n}^{-1}=q^{-a_{nm}}F_{m}\,,
\ee
\ese
\be\label{slnc} E_{n}F_{m}-F_{m}E_{n}=\delta_{nm}
\frac{K_n -K^{-1}_n}{q-q^{-1}}\,,
\ee
\be\label{slnd}
\sum\limits_{r=0}^{1-a_{nm}}(-1)^{r}
\left[\begin{array}{c}
1-a_{nm} \\
r
\end{array}\right]_{q}E_{n}^{1-a_{nm}-r} E_{m}E_{n}^{r}=0 ,\ \  \mbox{if}\;
n\neq
m,\\
\sum\limits_{r=0}^{1-a_{nm}}(-1)^{r}
\left[\begin{array}{c}
1-a_{nm} \\
r
\end{array}\right]_{q}F_{n}^{1-a_{nm}-r} F_{m}F_{n}^{r}=0 , \ \ \mbox{if}\;
n\neq m.\\ \ee
Here we have used the standard notations
\be
\left[\begin{array}{c}
m \\ n
\end{array}\right]_{q}=\frac{[m]_{q}!}{[n]_{q}![m-n]_{q}!}\;,\ \ \ \
[m]_{q}!=\prod\limits_{1\leq j \leq m}\frac{q^{j}-q^{-j}}{q-q^{-1}}\;.
\ee
The largest $\CC (q)$--algebra, the simply-connected rational form
$U^{P}_{q}(\frak{sl}(N)),$ is obtained by adjoining to
$U_{q}^{Q}(\frak{sl}(N))$ the invertible elements $L_n
\,,n=1\ldots,N-1,$ such that $K_n =\prod_m L_{m}^{a_{mn}}.$ The
relations (\ref{slnb}) are  now replaced by
\be\label{scrp}
L_{n}E_{m}L_{n}^{-1}=q^{\delta_{nm}}E_{m}\,,\ \ \ \
L_{n}F_{m}L_{n}^{-1}=q^{-\delta_{nm}}F_{m}\,.
\ee
Denote by $M$  any lattice such that $Q \subseteq M\subseteq P$.
 A general
rational form  $U_{q}^{M}({\frak g})$ is obtained by
adjoining to $U^{Q}_{q}({\frak g})$ the elements
$K_{\beta_i}=\prod_{j}L_{j}^{m_{ij}} $ for every
 $\beta_i =\sum_{j}m_{ij}\lambda_{j}\in M$, where $\lambda_j$
are the fundamental weights. Note that various
 rational  forms $U_{q}^{M}$ have interesting applications to the
construction of quantum integrable systems of Toda type.


\subsection{$U_q(\frak{g})$ in terms of the quantum tori}

We start with the definition of quantum torus algebra
$\hat\TT_q$. Let $\hat\TT_q$ be the associative $\CC
(q)$--
algebra of the rational functions of invertible elements ${\bf v}_{nj},\,
n=1,\ldots,N;\,j=1,\ldots,n$ and
${\bf u}_{nj}$, $n=1,\ldots,N-1;\,j=1,\ldots,n$ which are subjected to the
relations
\be\label{tq1}
{\bf v}_{nj}{\bf v}_{mk}={\bf v}_{mk}{\bf v}_{nj}\,,\ \ \ \
{\bf u}_{nj}{\bf u}_{mk}={\bf u}_{mk}{\bf u}_{nj}\,,\\
{\bf u}_{nj}{\bf v}_{mk}=q^{\delta_{nm}\delta_{jk}}{\bf v}_{mk}
{\bf u}_{nj}\,.
\ee
\begin{te}\label{gkth}
(i) Let $K_{nn}^{\pm 1}$, $E_{nm}$ be the generators of \ $U_q(\frak{gl}(N))$.
The following explicit expressions define the embedding $\pi$:
$U_q (\frak{gl}(N))\hookrightarrow\hat\TT_{q}$
\be\label{nc1}
\pi(K_{nn})=\prod\limits_{j=1}^{n}{\bf v}_{nj}\prod\limits_{j=1}^{n-1}
{\bf v}^{-1}_{n-1,j}\,,\\
\hspace{-0.5cm}
\pi(E_{n,n+1})=\!-\frac{q^{-1}}{q-q^{-1}}\prod\limits_{j=1}^{n+1}
{\bf v}_{n+1,j}^{-1}\prod\limits_{j=1}^{n}{\bf v}_{n,j}
\sum_{j=1}^n{\bf v}_{nj}^{-3}\,\frac{\prod\limits_{r=1}^{n+1}
({\bf v}_{nj}^2 -
q{\bf v}_{n+1,r}^2)} {\prod\limits_{s\neq j}({\bf v}_{nj}^2-
{\bf v}_{ns}^2)}\,{\bf u}_{nj}^{-1}\,,\\
\hspace{-0.5cm}
\pi(E_{n+1,n})=\frac{1}{q-q^{-1}}\prod\limits_{j=1}^{n}
{\bf v}_{nj}\prod\limits_{j=1}^{n-1}{\bf v}^{-1}_{n-1,j}\,
\sum_{j=1}^n {\bf v}_{nj}^{-1}
\frac{\prod\limits_{r=1}^{n-1}({\bf v}_{nj}^2 -
q^{-1}{\bf v}_{n-1,r}^2)}
{\prod\limits_{s\neq j}({\bf v}_{nj}^2-{\bf v}_{ns}^2)}\,
{\bf u}_{nj}\,.
\ee

(ii) Let $L^{\pm1}_n,K^{\pm1}_n,E_{n},F_{n}$ be the generators of
$U_q(\frak{sl}(N))$. The following explicit expressions define the
embedding  $\pi$: $U_q(\frak{sl}(N)) \hookrightarrow
\hat\TT_{q}$:
\be\label{rsln1}
\pi(L_{n})=\prod\limits_{j=1}^{n}{\bf v}_{nj}\,,\ \ \ \
\prod\limits_{j=1}^{N}{\bf v}_{\N j}=1\,,\\
\pi(K_{n})=\prod\limits_{m}(\prod\limits_{j=1}^{m}{\bf
v}_{mj})^{a_{mn}}\;\;,
\ee
\be\label{rsln3}
\hspace{-0.3cm}
\pi(E_{n})=-\frac{q^{-1}}{q-q^{-1}}\prod\limits_{j=1}^{n+1}
{\bf v}_{n+1,j}^{-1}\prod\limits_{j=1}^{n}{\bf v}_{nj}\sum_{j=1}^n
{\bf v}_{nj}^{-3}\,\frac{\prod\limits_{r=1}^{n+1}
({\bf v}_{nj}^2-q{\bf v}_{n+1,r}^2)}
{\prod\limits_{s\neq j}({\bf v}_{nj}^2 -{\bf v}_{ns}^2)}\,
{\bf u}_{nj}^{-1}\,,
\ee
\be\label{rsln4}
\pi(F_{n})=\frac{1}{q-q^{-1}}\prod\limits_{j=1}^{n}
{\bf v}_{nj}\prod\limits_{j=1}^{n-1}{\bf v}^{-1}_{n-1,j}
\sum_{j=1}^n {\bf v}_{nj}^{-1}\,\frac{\prod\limits_{r=1}^{n-1}
({\bf v}_{nj}^2-q^{-1}{\bf v}_{n-1,r}^2)}
{\prod\limits_{s\neq j}({\bf v}_{nj}^2 -{\bf v}_{ns}^2)}\,
{\bf u}_{nj}\,.
\ee
\end{te}
The algebra $\hat\TT_{q}$  is very close to the minimal
skew field of fractions of $U_q(\frak{gl}(N))$  for which the
inclusion of the universal enveloping algebra is possible.
Consider the skew field of fractions ${\cal D}(U_{q}({\frak g}))$
of $U_{q}({\frak g})$ which consists of the elements of the form
$u\cdot v^{-1}$ or $x^{-1}\cdot y$, where
$u,v,x,y\in U_q({\frak g})$ (for more details see e.g. \cite{GK}).
It appears that the skew field of fractions ${\cal D}( U_{q}(\frak{gl}(N)))$
is obtained by partial symmetrization of  the algebra
$\hat\TT_{q}$
\be
{\cal D}(U_{q}(\frak{gl}(N)))=(\hat\TT_q)^{{\otimes}_{n=1}^NS_{n}}
\ee
where, for $n<N$, the group $S_n$ acts as:
\bse\label{sim1}
\be\label{sim1a}
\sigma\!:\, {\bf v}_{nj}\rightarrow {\bf v}_{n,\sigma(j)},
\ee
\be\label{sim1b}
\sigma\!:\, {\bf u}_{nj}\rightarrow {\bf u}_{n,\sigma(j)},
\ee
\ese
and the group $S_\N$ acts as:
\be\label{sim2}
\sigma\!:\, {\bf v}_{\N j}\rightarrow {\bf v}_{\N,\sigma(j)}.
\ee
This proposition can be considered as a generalization of the
Gelfand-Kirillov  theorem for $U(\mathfrak{g})$ \cite{GK}.

\subsection{ $U_{q}(\frak{g})\,\otimes U_{\tilde q}(
\frak{\check{g}})$-bimodule structure}

It turns out that  the  consideration of  the infinite-dimensional
representation of various  rational forms
 (\ref{sln})-(\ref{scrp}) of quantum
groups reveals  new phenomena. Namely some representations of the
quantum group $U_q(\frak{g})$ poses a natural structure of a
$U_{q}(\frak g)\,\otimes U_{\tilde q} (\frak {\check{g}})$-bimodule
where $\check {\mathfrak{g}}$ is Langlands dual to $\mathfrak{g}$ and
$\log q/2\pi i=-(\log{\tilde q}/2\pi i)^{-1}$.
The key point is the interpretation of the appropriately defined centralizer
of the image of $U_q(\mathfrak{g})$ in these representations
in terms of the representation of $U_{\tilde q} (\frak{\check{g}})$.
The presence of the Langlands dual will be verified below in the case of
$U_q(\mathfrak{sl}(N))$. It turns out that the centralizer construction
leads to the connection between various forms of $U_q(\mathfrak{sl}(N))$.
In particular, when we start from the minimal (adjoint) form,  we come to \
what may be called the maximal form over the quadratic extension,
which is in agreement with the classical
picture of the duality between $SL(N)$ and $PSL(N)$.

Let us begin  with the construction of a representation of
$U_q(\mathfrak{gl}(N))$  generalizing the representation of
$U(\mathfrak{gl}(N))$ described in Section 2. Let $\mathcal{V}$ be
the space of functions of $\gamma_{nj}$ with $n=1,\ldots,N-1;\,1\leq j\leq n$
and $\mathcal{V}^s$ be the space of functions of $\gamma_{nj}$ invariant
under the action of $\otimes_{n=1}^{N-1}S_n$ according to (\ref{sim1a}).
Let $\mathcal{R}$ be the algebra of rational functions of the exponents of
the linear functions of $\gamma_{nj}$ with $n=1,\ldots N;\,1\leq j\leq n$
and $\d_{\gamma_{nj}}:=\frac{\pr}{\pr \gamma_{nj}}$ with
$n=1,\ldots N-1;\,1\leq j\leq n$ and $\mathcal{R}^s$ be the algebra of
rational functions of the exponents of the linear functions of
$\gamma_{nj}$ and $\d_{\gamma_{nj}}$
invariant under the action of $\otimes_{n=1}^NS_n$ according to
(\ref{sim1}), (\ref{sim2}).
One has the following embedding
$\varrho:\hat\TT_{q}\hookrightarrow \mathcal{R}$
\be\label{tor1}
\varrho({\bf u}_{nj})=e^{i\omega_1 \d_{\gamma_{nj}}}\,,
\ \ \ \varrho({\bf v}_{nj})=e^{\frac{2\pi\gamma_{nj}}{\omega_2}}\,,
\ \ \ \varrho(q)=e^{\frac{2\pi i\omega_1}{\omega_2}}.
\ee
Using Theorem \ref{gkth} we obtain the representation of
$U_{q}(\frak{gl}(N))$ in terms of the difference operators.
By an abuse of  notation we shall also  denote by $\varrho$ the composition
$\varrho\circ\pi$.
\begin{prop}\label{langl}
The following expressions define a representation
of $U_q(\frak{gl}(N))$ with $q=e^{\frac{2\pi
i\omega_1}{\omega_2}}$:
\be\label{wnc1}
\varrho(K_{nn})=e^{\frac{2\pi}{\omega_2}\,\big(\sum\limits_{j=1}^n\gamma_{nj}-
\sum\limits_{j=1}^{n-1}\gamma_{n-1,j}\big)}\,,\\
\varrho(E_{n,n+1})=\frac{2ie^{\frac{\pi i\omega_1}{\omega_2}(n-1)}}
{\sin\frac{\textstyle 2\pi\omega_1}{\textstyle\omega_2}}
\sum_{j=1}^n\frac{\prod\limits_{r=1}^{n+1}
\sinh\frac{2\pi}{\omega_2}(\gamma_{nj}-\gamma_{n+1,r}-
\frac{i\omega_1}{2})}
{\prod\limits_{s\neq j}\sinh\frac{2\pi}{\omega_2}(\gamma_{nj}-\gamma_{ns})}
\,e^{-i\omega_1\d_{\gamma_{nj}}}\,,\\
\hspace{-0.8cm}
\varrho(E_{n+1,n})=-\frac{i e^{-\frac{\pi i \omega_1}{\omega_2}(n-1)}}
{2\sin\frac{\textstyle 2\pi\omega_1}{\textstyle\omega_2}}
\sum_{j=1}^n\frac{\prod\limits_{r=1}^{n-1}
\sinh\frac{2\pi}{\omega_2}(\gamma_{nj}-\gamma_{n-1,r}+
\frac{i\omega_1}{2})}
{\prod\limits_{s\neq j}\sinh\frac{2\pi}{\omega_2}(\gamma_{nj}-\gamma_{ns})}
\,e^{i\omega_1\d_{\gamma_{nj}}}\,.
\ee
\end{prop}
Consider the dual quantum torus $\hat\TT_{\tilde q}$ using the
invertible elements $\tilde{\bf v}_{nj},\, n=1,\ldots,N$; $j=1,\ldots,n$ and
$\tilde{\bf u}_{nj},\,n=1,\ldots,N-1;\,j=1,\ldots,n$ subjected to the
relations
\be
\tilde{\bf v}_{nj}\tilde{\bf v}_{mk}=\tilde{\bf v}_{mk}\tilde{\bf v}_{nj}\,,
\ \ \ \
\tilde{\bf u}_{nj}\tilde{\bf u}_{mk}=
\tilde{\bf u}_{mk}\tilde{\bf u}_{nj}\,,\\
\tilde{\bf u}_{nj}\tilde{\bf v}_{mk}=\tilde q^{\;\delta_{nm}\delta_{jk}}
\tilde{\bf v}_{mk}\tilde{\bf u}_{nj}\,.
\ee
One has the dual embedding
$\tilde\varrho:\hat\TT_{\tilde q}\hookrightarrow
\mathcal{R}$
\be\label{tor2}
\tilde\varrho(\tilde{\bf u}_{nj})=e^{i\omega_2 \d_{\gamma_{nj}}}\,,
\ \ \ \tilde\varrho(\tilde{\bf v}_{nj})=
e^{-\frac{2\pi\gamma_{nj}}{\omega_1}}\,,\ \ \
\tilde\varrho(\tilde q)=e^{-\frac{2\pi i\omega_2}{\omega_1}}.
\ee
and the two actions of the torus and its dual  on the same space of
functions of $\gamma_{ni}$ mutually commute. Hence, we have the following
\begin{prop}
The operators
\be\label{dg2}
\tilde\varrho(K_{nn})=e^{-\frac{2\pi }{\omega_1}\,
\big(\sum\limits_{j=1}^n\gamma_{nj}-
\sum\limits_{j=1}^{n-1}\gamma_{n-1,j}\big)}\,,\\
\tilde\varrho(E_{n,n+1})=-\frac{2i\,e^{-\frac{\pi i\omega_2}{\omega_1}(n-1)}}
{\sin\frac{\textstyle 2\pi\omega_2}{\textstyle\omega_1}}
\sum_{j=1}^n\frac{\prod\limits_{r=1}^{n+1}
\sinh\frac{2\pi}{\omega_1}(\gamma_{nj}-\gamma_{n+1,r}-
\frac{i\omega_2}{2})}
{\prod\limits_{s\neq j}\sinh\frac{2\pi}{\omega_1}(\gamma_{nj}-\gamma_{ns})}
\,e^{-i\omega_2\d_{\gamma_{nj}}}\,,\\
\tilde\varrho(E_{n+1,n})=\frac{i\,e^{\frac{\pi i\omega_2}{\omega_1}(n-1)}}
{2\sin\frac{\textstyle 2\pi\omega_2}{\textstyle\omega_1}}
\sum_{j=1}^n\frac{\prod\limits_{r=1}^{n-1}
\sinh\frac{2\pi}{\omega_1}(\gamma_{nj}-\gamma_{n-1,r}+
\frac{i\omega_2}{2})}{\prod\limits_{s\neq j}
\sinh\frac{2\pi}{\omega_1}(\gamma_{nj}-\gamma_{ns})}
\,e^{i\omega_2\d_{\gamma_{nj}}}\,,
\ee
generate a representation of the dual quantum group
$U_{\tilde{q}}(\frak{gl}(N))$, where
${\tilde q}=e^{-\frac{2\pi i \omega_2}{\omega_1}}$. The operators
$\varrho(U_{q}(\frak{gl}(N)))$ and $\tilde\varrho(U_{\tilde q}(\frak{gl}(N)))$
are commute by construction.
\end{prop}
Thus we have a structure of  a $U_{q}(\frak{gl}(N))\otimes
U_{\tilde{q}}(\frak{gl}(N))$-bimodule. More precisely, this bimodule
may be characterized by the condition
\be\label{constr}
\hspace{2cm}
\varrho(K_{nn})=\tilde\varrho({K}_{nn})^{-\tau}\,,\ \ \ \
n=1,\ldots,N\,,
\ee
where $\tau=\omega_1/\omega_2$. We will discuss the better way to formulate
this condition latter in this section.

Let us remark that we actually constructed the embeddings
$\varrho:U_q(\frak{gl}(N))\hookrightarrow \mathcal{R}^s$
and $\tilde{\varrho}:U_{\tilde q}
(\frak{gl}(N))\hookrightarrow \mathcal{R}^s$. Consider the
centralizer $[{\cal D}( \varrho(U_{q}\frak{gl}(N)))]'$ of the
algebras ${\cal D}( \varrho(U_{q}\frak{gl}(N)))$ as a subalgebra
of $\mathcal{R}^s$.  Then we have
\begin{te}\label{commt}
\be
[{\cal D}( \varrho(U_{q}\frak{gl}(N)))]'= {\cal D}(\tilde\varrho
(U_{\tilde q} \frak{gl}(N))). \ee
\end{te}
This Theorem  is proved by direct calculation.

Consider now more interesting case of $U_{q}(\frak{sl}(N))$. It appears
that the centralizer of the minimal (adjoint) rational form
of the quantum group $U_q(\mathfrak{sl}(N))$ is
described in terms of the different forms of the same quantum group.
  This may be considered as an
indication of the fact that the Langlands dual of $PSL(N)$ is
$SL(N)$. One may conjecture that the  Langlands dual quantum group
$U_{\tilde q}(\check{{\frak g}})$  may be generally obtained  as a  result
of explicit calculations of the centralizer in the appropriate
generalization of the above construction to arbitrary quantum
groups.

Below we give the explicit formulas for
$U_q(\frak{sl}(2))$. The adjoint rational form of the
quantum group $U_q^{Q}(\frak{sl}(2))$ is
generated by elements $K,K^{-1},E,F$ subjected to the  relations
\be\label{ad}
K E K^{-1}=q^{2}E\,,\ \ \ K F K^{-1}=q^{-2} F\,,\\
EF -F E=\frac{K-K^{-1}}
{q-q^{-1}}\,.
\ee
The representation  is given by
\be
\varrho(K)=e^{\frac{4\pi\gamma_{11}}{\omega_2}}
\,,\\
\hspace{-0.5cm}
\varrho(E)=\frac{2i}
{\sin\frac{\textstyle 2\pi\omega_1}{\textstyle\omega_2}}
\sinh\frac{2\pi}{\omega_2}(\gamma_{11}-\nu-\frac{i\omega_1}{2})
\sinh\frac{2\pi}{\omega_2}(\gamma_{11}+\nu-\frac{i\omega_1}{2})\,
e^{-i\omega_1 \d\gamma_{11}},\\
\varrho(F)=-\frac{i}{\sin\frac{\textstyle 2\pi\omega_1}{\textstyle\omega_2}}
\,e^{i\omega_1 \d\gamma_{11}}\,,\\
\varrho(q)=e^{\frac{2\pi i\omega_1}{\omega_2}}\,.
\ee
This representation is obtained from those of
$U_{q}(\frak{gl}(2))$ by the restriction $\gamma_{21}=-\gamma_{22}:=\nu$.
The centralizer in $\cal R$ is generated by the algebra of functions
of the dual torus $\TT_{\tilde q^{\,1/2}}$:
\be
{\tilde u}{\tilde v}=\tilde q^{\,1/2}{\tilde v}{\tilde u}\,\,,
\ee
where
\be
\tilde\varrho(\tilde u)=e^{\frac{i\omega_2}{2}\d\gamma_{11}},\ \ \ \
\tilde\varrho(\tilde v)=e^{-\frac{2\pi\gamma_{11}}{\omega_1}},\ \ \ \ \
\tilde\varrho(\tilde  q)=e^{-\frac{2\pi i\omega_2}{\omega_1}}.
\ee
The  centralizer in ${\cal R}^s$ may be described as an image
of the skew field of
fractions of the following algebra generated by
${\tilde L},{\tilde E}, {\tilde F}, {\tilde K}={\tilde L}^4$
over the quadratic extension $\CC ({\tilde q}^{1/2})$:
\be
\tilde L\tilde E{\tilde L}^{-1}=\tilde q^{\,1/2}\tilde E\,,\ \ \ \
\tilde L\tilde F{\tilde L}^{-1}=\tilde q^{\,-1/2}\tilde F , \\
\tilde E\tilde F-\tilde F\tilde E=
\frac{ \tilde K-{\tilde K}^{-1}}{\tilde q-\tilde q^{-1}}\, .
\ee
under the representation
\be
\tilde\varrho(L)=e^{-\frac{2\pi\gamma_{11}}{\omega_1}},\\
\hspace{-0.6cm}
\tilde\varrho(E)=-\frac{2i}
{\sin\frac{\textstyle2\pi\omega_2}{\textstyle\omega_1}}
\sinh\frac{2\pi}{\omega_1}(2\gamma_{11}-\nu-\frac{i\omega_2}{2})
\sinh\frac{2\pi}{\omega_1}(2\gamma_{11}+\nu-\frac{i\omega_2}{2})\,
e^{-\frac{i\omega_2}{2}\d\gamma_{11}},\\
\tilde\varrho(F)=\frac{i}
{2\sin\frac{\textstyle 2\pi\omega_2}{\textstyle\omega_1}}\,
e^{\frac{i\omega_2}{2}
\d\gamma_{11}}.
\ee
This algebra may be considered as a  maximal form of
$U_{\tilde q}(\mathfrak{sl}(2))$ over the quadratic extension
$\CC ({\tilde q}^{1/2})$.
Let us stress that  the reconstruction of the algebra from its skew field of
fractions is not unique. Thus, in this example the same centralizer may
be interpreted as the image of the skew field of fractions of the
simply-connected form $U^P_{\tilde q^{\,1/2}}(\frak{sl}(2))$:
\be
\tilde\varrho(L)=e^{-\frac{2\pi\gamma_{11}}{\omega_1}},\\
\hspace{-0.6cm}
\tilde\varrho(E)=-\frac{2i}
{\sin\frac{\textstyle\pi\omega_2}{\textstyle\omega_1}}
\sinh\frac{2\pi}{\omega_1}(\gamma_{11}-\nu-\frac{i\omega_2}{4})
\sinh\frac{2\pi}{\omega_1}(\gamma_{11}+\nu-\frac{i\omega_2}{4})\,
e^{-\frac{i\omega_2}{2}\d\gamma_{11}},\\
\tilde\varrho(F)=\frac{i}
{2\sin\frac{\textstyle\pi\omega_2}{\textstyle\omega_1}}\,
e^{\frac{i\omega_2}{2}
\d\gamma_{11}}.
\ee
where the operators satisfy the  relations:
\be
\tilde L\tilde E {\tilde L}^{-1}=\tilde q^{\,1/2}\tilde E\,,\ \ \ \
\tilde L \tilde F {\tilde L}^{-1}=\tilde q^{\,-1/2}\tilde F , \\
\tilde E\tilde F -\tilde F\tilde E=
\frac{ \tilde K -{\tilde K}^{-1}}{\tilde q^{\,1/2}-\tilde q^{\,-1/2}}\,,
\ee
and $\tilde K={\tilde L}^2$.  Moreover, there is an
isomorphism of the algebras $\mathcal{D}(U^P_{\tilde q{\,1/2}}(\frak{sl}(2)))$ and
$\mathcal{D}(U^Q_{\tilde q{\,1/4}}(\frak{sl}(2)))$. This allows to
reformulate the duality in a more symmetric form. Taking $p=e^{\pi i\tau}$ and $\tilde
p=e^{-\pi i/\tau}$ with  $\tau=\frac{2\omega_1}{\omega_2}$ one gets
the duality  between  algebra $U^Q_p(\frak{sl}(2))$ and
$U^Q_{\tilde p}(\frak{sl}(2))$ which leads to the modular double
considered in \cite{Fad}, \cite{Fad1}. Let us stress however that the
dual quantum deformation
parameters enter here in a non-standard way.

Finally, consider the algebra $U^P_q(\frak{sl}(2))$ such that
${\cal D}(U^P_q(\frak{sl}(2)))=
\TT_q:=\{v=e^{\frac{2\pi\gamma_{11}}{\omega_2}},\,
u=e^{i\omega_1\d_{\gamma_{11}}}\}$ with
$q=e^{\frac{2\pi i\omega_1}{\omega_2}}$. Obviously,
$[{\cal D}(U^P_q(\frak{sl}(2)))]'=\TT_{\tilde q}:=
\{\tilde v=e^{-\frac{2\pi\gamma_{11}}{\omega_1}},\,
\tilde u=e^{i\omega_2\d_{\gamma_{11}}}\}$ with
$\tilde q=e^{-\frac{2\pi i\omega_2}{\omega_1}}$. In other words, the
algebras $\mathcal{D}(U^P_q(\frak{sl}(2)))$ and
$\mathcal{D}(U^P_{\tilde q}(\frak{sl}(2)))$
with $q=e^{2\pi i\tau}$ and $\tilde q=e^{-2\pi i/\tau}$,
$\tau=\frac{\omega_1}{\omega_2}$
centralize each other.

Clearly to use the centralizers as a way to describe the Langlands
dual pairs deserves additional structures on the representation space
comparing to what was discussed above. We are going to consider  these
matters elsewhere.
Let us also remark that the use of the continuous powers of the Cartan
generators in (\ref{constr}) is not quite appropriate in our
setting. A possible  solution is to identify instead the actions
of the centers of both algebras. It can be shown that the center
is described in terms of the symmetric polynomials of
${\bf v}_{\N i}$ and thus  the bimodule structure may be described
equivalently as
\be\label{constr31}
\varrho({\bf v}_{\N j})={\tilde\varrho}({\tilde{\bf v}}_{\N j})^{-\tau}\,.
\ee
If we consider $\log{\bf v}_{\N j}$ as legitimate operators then
the relations (\ref{constr31}) make sense. However we believe that
a proper description  of the structure of this bimodule which
solves this problem should be given in a different way. Let us
notice that the description of the universal enveloping algebra
in terms of  the quantum tori has obvious asymmetry. The variables
${\bf v}_{Ni}$ are coordinate functions on the commutative sub-torus
and thus there is no natural definition of the dual torus through
the centralizer. Thus it is natural to guess  that some
generalization of the universal enveloping algebra  provides
 the proper setting for the discussion of the Langlands duality for
quantum groups through the centralizers. The most natural
candidate is the quantum group analogs $\mathcal{A}(G_q)$ of the
differential operators on the group $Diff(G)$ and on the basic
affine space $Diff(G/N)$. This leads to the  interpretation of the
resulting $\mathcal{A}(G_q)\otimes
\mathcal{A}(\check{G}_{\tilde{q}})$-bimodule as an explicit
realization of the Morita equivalence of the algebras
$\mathcal{A}(G_q)$ and $\mathcal{A}(\check{G}_{\tilde{q}})$.
 Preliminary results based on the generalized Gelfand-Zetlin
representation  (see \cite{GKL2}) support this conjecture. We are
going to discuss this approach in the future.

\subsection{Whittaker modules}

Let us define the special class of the representations of
$U_{q}(\frak{gl}(N))\otimes U_{\tilde{q}}(\frak{gl}(N)),$
generalizing the Whittaker module for the classical algebras
described in the previous section. The Whittaker vectors for
$U_{q}(\frak{gl}(N))\otimes U_{\tilde q}(\frak{gl}(N))$-bimodule are
defined using the generalization of the definition in \cite{KLS}.
\footnote{
For the case of $U_{q}(\frak g )$, where $\frak g$ is an arbitrary
simple Lie algebra, the construction of the non-degenerate
characters of nilpotent subalgebras was done by Sevostyanov \cite{Sev}.
However, the appearance of the bimodule structure
reveals larger symmetry in the representation theory of quantum groups.}
\begin{de}
The Whittaker vectors $w_{\N}$ and $w'_{\N}$ are defined by equations
\be\label{wv1}
E_{n,n+1}w_\N=\frac{\chi_n}{q-q^{-1}}\;
K_{nn}\prod_{m=1}^NK_{mm}^{c_{nm}-c_{n+1,m}}w_\N\,,\\
\tilde E_{n,n+1}w_\N=\frac{\chi_n}{\tilde q-\tilde q^{-1}}\;
\tilde K_{nn}\prod_{m=1}^N\tilde K_{mm}^{c_{nm}-c_{n+1,m}}w_\N\,,
\ee
\be\label{wv2}
E_{n+1,n}w'_\N=\frac{\chi'_n}{q-q^{-1}}\;
K_{nn}^{-1}\prod_{m=1}^NK_{mm}^{c'_{nm}-c'_{n+1,m}}w'_\N\,,\\
\tilde E_{n+1,n}w'_\N=\frac{\chi'_n}{\tilde q-\tilde q^{-1}}\;
\tilde K_{nn}^{-1}\prod_{m=1}^N\tilde K_{mm}^{c'_{nm}-c'_{n+1,m}}w'_\N\,,
\ee
where $||c_{nm}||,\,||c_{nm}'||$ are the $N\times N$ symmetric matrices
such that $c_{nm}-c_{n+1,m}$, $c'_{nm}-c'_{n+1,m}$ are integers, and
$\chi_n,\,\chi_n'$ are arbitrary parameters.
\end{de}
The direct calculations show that the following statement is true:
\begin{prop}
The defining equations (\ref{wv1}), (\ref{wv2}) are
consistent with the full set of the Serre relations (\ref{d3}) and their
dual analogues.
\end{prop}

The results of Section 2 can be naturally extended to the quantum
group case. The structure of the Whittaker vectors and Whittaker
modules remains essentially the  same. In particular, the
Whittaker vectors can be written in a form similar to equations
(\ref{wv}), (\ref{wv'}). For example, there is a solution to
(\ref{wv2}) which is unique up to multiplication
by an arbitrary double-periodic function:
\be\label{lw3}
w'_\N=e^{-\frac{\pi i}{\omega_1\omega_2}\sum\limits_{n,m=1}^Nc'_{nm}h_nh_m}
\prod_{n=1}^{N-1} e^{\frac{\pi
i}{\omega_1\omega_2}\sum\limits_{p=1}^n\gamma_{np}^2
+\frac{\pi(\omega_1+\omega_2)d_n}
{\omega_1\omega_2}\sum\limits_{p=1}^n\gamma_{np}}\,,
\ee
where
\be
h_n=\sum_{j=1}^n\gamma_{nj}-\sum_{j=1}^{n-1}\gamma_{n-1,j}\,,
\ee
and $d_n=2n-c'_{nn}+2c'_{n,n+1}-c'_{n+1,n+1}-1$. In the present example
$\chi'_n=(-1)^{d_n}$.

Let us denote by  $W'_\N$  the Whittaker module  with the
cyclic vector (\ref{lw3}). It can be proved that it is spanned by the
symmetric polynomials of the variables
$e^{\pm\frac{2\pi\gamma_{nj}}{\omega_1}},\,
e^{\pm\frac{2\pi\gamma_{nj}}{\omega_2}}$ (compare with Section 2.2).
The matrix elements of the particular group elements between
the Whittaker vectors leads to the explicit expressions for the
wave functions of the generalized quantum Toda theories and will
considered elsewhere.

\section{The QISM's eigenfunction via representation theory}
\subsection{An infinite dimensional representations of the  Yangian}
The aim of this section  is to introduce a special type of an infinite
dimensional  representations  of the $Y(\frak{gl}(N))$.
 This allows to connect the QISM methods of the solution of the
 integrable  system based on the representation theory of Yangian
and the solution  based on the representation discussed in Section
2.

We start with  some well known facts of the Yangian theory
\cite{Dr2}, \cite{Dr1} (see also the recent review \cite{MNO}).
The Yangian $Y(\frak{gl}(N))$ is an associative Hopf algebra
generated by the elements $T_{ij}^{(r)}$, where $i,j=1,\ldots,N$
and $r=0,1,2,\ldots$, subject the following relations. Consider
the $N\times N$ matrix $T(\la)=||T_{ij}(\la)||_{i,j=1}^\N$ with
operator-valued entries
\be\label{y1}
T_{ij}(\la)=\la\delta_{ij}+\sum_{r=0}^\infty T_{ij}^{(r)}\la^{-r}.
\ee
Let
\be\label{r-mat}
\hspace{2cm}
R(\la)\,=\,I\otimes I+i\hbar P/\la\;,\hspace{1cm}
P_{ik,jl}=\delta_{il}\delta_{kj}\,,
\ee
be an $N^2\times N^2$ numerical matrix (the Yang $R$-matrix). Then
the relations between the generators $T_{ij}^{(r)}$ can be written
in the standard form
\be\label{y2}
R(\la-\mu)(T(\la)\otimes I)(I\otimes T(\mu))=
(I\otimes T(\mu))(T(\la)\otimes I)R(\la-\mu)\,.
\ee
The centre of the Yangian is generated by the
coefficients of the formal Laurent series (the quantum
determinant of $T(\la)$ in the sense of \cite{KS}):
\be\label{y4}
{\rm det}_q T(\la)\\ =\,\sum_{s\in S_\N}{\rm sign}\,s\;
T_{s(1),1}(\la-i\hbar\rho^{(\N)}_1)\ldots
T_{s(k),k}(\la-i\hbar\rho^{(\N)}_k) \ldots
T_{s(\N),\N}(\la-i\hbar\rho^{(\N)}_\N)\,, \hspace{-0.5cm}
\ee
where $\rho^{(\N)}_n=\ts{1}{2}(N-2n+1),\,n=1,\ldots,N$ and the
summation is  over the elements of the permutation group $S_\N$.
Let $X(\la)=||X_{ij}(\la)||_{i,j=1}^n$ be an $n\times n$
submatrix of the matrix $||T_{ij}(\la)||_{i,j=1}^\N$. It is
obvious from the explicit form of $R_\N(\la)$ that this submatrix
satisfies an analogue of relations (\ref{y2}). The quantum
determinant ${\rm det}_q X(\la)$ is defined similarly to
(\ref{y4}) (with the evident change $N\to n$).

The following way to describe the Yangian $Y(\frak{gl}(N))$ was
introduced in \cite{Dr2}. Let ${\bf A}_n(\la)$, $n=1,\ldots, N$,
be the quantum determinants of the submatrices, determined by the
first $n$ rows and columns, and let  the operators ${\bf
B}_n(\la),{\bf  C}_n(\la)$, $n=1,\ldots, N-1$,  be the quantum
determinants of the  submatrices with  elements $T_{ij}(\la)$,
where $i=1,\ldots,n$; $j=1,\ldots,n-1,n+1$ and
$i=1,\ldots,n-1,n+1$; $j=1,\ldots,n$, respectively. The expansion
coefficients of $ {\bf A}_n(\la),{\bf B}_n(\la) ,{\bf C}_n(\la)$,
$n=1,\ldots,N-1$, with respect to $\la$, together with those of
${\bf A}_\N(\la)$, generate the algebra $ Y(\frak{gl}(N)) $. The
parts of the relations, which we use below, are as follows:
\be\label{cw1}
\hspace{3cm} [{\bf A}_n(\la),{\bf A}_m(\mu)]=0\;;
\hspace{1cm}(n,m=1,\ldots, N),\\
\hspace{2cm}[{\bf B}_n(\la),{\bf B}_m(\mu)]=0\,;\;\ \
[{\bf C}_n(\la),{\bf C}_m(\mu)]=0\;;\hspace{1cm} (m\neq n\pm 1),\\
(\la-\mu+i\hbar){\bf A}_n(\la){\bf B}_n(\mu)\,=\,(\la-\mu)
{\bf B}_n(\mu){\bf A}_n(\la)+i\hbar{\bf A}_n(\mu){\bf B}_n(\la),\\
(\la-\mu+i\hbar){\bf A}_n(\mu){\bf C}_n(\la)\,=\, (\la-\mu)
{\bf C}_n(\la){\bf A}_n(\mu)+i\hbar{\bf A}_n(\la){\bf C}_n(\mu).
\ee

Let $A(\frak{gl}(N))$ be the commutative
subalgebra of $Y(\frak{gl}(N))$ generated by ${\bf
A}_n(\la),\,n=1,\ldots,N$. It was proved in \cite{Cher2} that
$A(\frak{gl}(N))$ is the maximal commutative subalgebra of
$Y(\frak{gl}(N))$.

There is  another realization of the  $Y(\frak{gl}(N))$ introduced
by Drinfeld  \cite{Dr1}.  The algebra $Y(\frak{gl}(N))$ is
generated by the coefficients of the formal series
\be
k_n(\la)=\la+\sum_{a=0}^\infty k_n^{(a)}\la^{-a}\,,\\
e_n(\la)=\sum_{a=0}^\infty e_n^{(a)}\la^{-a-1}\,,\\
f_n(\la)=\sum_{a=0}^\infty f_n^{(a)}\la^{-a-1}\,,\\
\ee
subjected to the commutation relations
\be\label{first}
[k_n(\la),k_m(\mu)]=0\,,\\

[k_n(\la),e_m(\mu)]\,=\,i\hbar
(\delta_{nm}-\delta_{n,m+1})\,k_n(\la)
\frac{e_m(\la)-e_m(\mu)}{\la\,-\,\mu} \,,\\

[k_n(\la),f_m(\mu)]\,=\,-\,i\hbar(\delta_{nm}-\delta_{n,m+1})
\frac{f_m(\la)-f_m(\mu)}{\la\,-\,\mu}\,k_n(\la)\,, \\

[e_n(\la),f_m(\mu)]=i\hbar
\frac{k_n^{-1}(\mu)k_{n+1}(\mu)-k_n^{-1}(\la)k_{n+1}(\la)}
{\la-\mu}\,\delta_{nm} \,,\\

[e_n^{(a+1)},e_m^{(b)}]-[e_n^{(a)},e_m^{(b+1)}]=
\ts{\imath\hbar}{2}\,
a_{nm}(e_n^{(a)}e_m^{(b)}+e_m^{(b)}e_n^{(a)})\,,\\

[f_n^{(a+1)},f_m^{(b)}]-[f_n^{(a)},f_m^{(b+1)}]=
-\ts{\imath\hbar}{2}\,
a_{nm}(f_n^{(a)}f_m^{(b)}+f_m^{(b)}f_n^{(a)})\,,\\

[e_n^{(a)},[e_n^{(b)},e_{n\pm 1}^{(c)}]]+
[e_n^{(b)},[e_n^{(b)},e_{n\pm 1}^{(c)}]]=0\,,\\

[f_n^{(a)},[f_n^{(b)},f_{n\pm 1}^{(c)}]]+
[f_n^{(b)},[f_n^{(b)},f_{n\pm 1}^{(c)}]]=0 \,,
\ee
where $a_{nm}=2\delta_{nm}-\delta_{n,m+1}-\delta_{n+1,n}$.  The
 relation between two realization is given by
\be
k_n(\la)=\frac{{\bf A}_n(\la-\ts{i(n-1)\hbar}{2})} {{\bf
A}_{n-1}(\la-\ts{in\hbar}{2})}\,,\\ e_n(\la)={\bf
A}^{-1}_n(\la-\ts{i(n-1)\hbar}{2}) {\bf
B}_n(\la-\ts{i(n-1)\hbar}{2})\,,\\ f_n(\la)={\bf
C}_n(\la-\ts{i(n-1)\hbar}{2}) {\bf
A}^{-1}_n(\la-\ts{i(n-1)\hbar}{2}) \,.\ee
Let us stress that Drinfeld's realization   is known for the
$Y(\frak g)$, where $\frak g$ be any simple Lie algebra
\cite{Dr2}-\cite{Dr1}.

There is a natural  epimorphism
${\pi}_\N :Y(\frak{gl}(N))\rightarrow{ U}(\frak{gl}(N))$
\be\label{bz1} \hspace{3cm}
\pi_{\N}(T_{jk}(\la))=\la\delta_{jk}-i\hbar E_{jk}\,,
\hspace{1cm}(j\,,\,k=1,\ldots, N).
\ee
Denote the images under $\pi_{\N}$ of the  generators ${\bf A}_n(\la)$ and
${\bf B}_n(\la),\,{\bf C}_n(\la)$  by $\ca_n(\la)$ and
$\cb_n(\la),\,\cc_n(\la)$,  respectively. Note that the images are the
polynomials in $\lambda$ of orders $n$ and $n-1$, respectively.

To obtain the representation of the Yangian $Y(\frak{gl}(N))$  we
start with the construction of a natural representation of the
Cartan subalgebra generated by $k_n(\lambda)$. It can be represented by the
rational functions as follows:
\be\label{az}
\prod\limits_{s=1}^{n}k_s(\lambda-i\hbar\rho^{(n)}_s)=
\prod\limits_{j=1}^{n}(\la-\gamma_{nj})\,,
\ee
where $\rho^{(n)}_s=\frac{1}{2}(n-2s+1)$ and $n=1,\ldots,N$.
Then the operators $A(\frak{gl}(N))$ act  by the polynomials
$\ca_{n}(\lambda)=\prod\limits_{j=1}^{n}(\la-\gamma_{nj}),\,
n=1,\ldots,N$.  We resolve  the rest of the  Yangian relations and
find the explicit expressions for the generators $e_n(\lambda)$, $
f_n(\lambda)$ (and $\cb_n(\la)$ and $\cc_n(\la)$) in terms of the
 operators acting on the space of functions depending on the
variables $\gamma_{nj}\,,j=1,\ldots, n;\,n=1,\ldots,N$.
\begin{te}\label{boz}
The operators
\be\label{new1}
k_{n}(\lambda)=
\frac{\prod\limits_{j=1}^n(\la-\gamma_{nj}-\ts{i(n-1)\hbar}{2})}
{\prod\limits_{j=1}^{n-1}(\la-\gamma_{n-1,j}-\ts{i n\hbar}{2})}\,,\\
e_{n}(\lambda)=
\sum_{j=1}^{n}\frac{1}{\lambda-\gamma_{nj}-\ts{i(n-1)\hbar}{2}}\;
\frac{\prod\limits_{r=1}^{n+1}(\gamma_{nj}-\gamma_{n+1,r}-\frac{i\hbar}{2})}
{\prod\limits_{s\neq j}(\gamma_{nj}-\gamma_{ns})}\,,\\
f_{n}(\lambda)=
\sum_{j=1}^{n}\frac{1}{\lambda-\gamma_{nj}-i\hbar-\ts{i(n-1)\hbar}{2}}\;
\frac{\prod\limits_{r=1}^{n-1}(\gamma_{nj}-\gamma_{n-1,r}+\frac{i\hbar}{2})}
{\prod\limits_{s\neq j}(\gamma_{nj}-\gamma_{ns})}\,
\,e^{i\hbar\d_{\gamma_{nj}}},
\ee
satisfy the complete set of relations (\ref{first}) and, therefore, define a
representation of the Yangian $Y(\frak{gl}(N))$.
\end{te}

As a consequence we have \be\label{bza}
\ca_{n}(\lambda)=\prod\limits_{j=1}^{n}(\la-\gamma_{nj}),\\
\cb_{n}(\la)=\sum_{j=1}^{n}\prod_{s\neq j}
\frac{\la-\gamma_{ns}}{\gamma_{nj}-\gamma_{ns}} \,
\prod_{r=1}^{n+1}(\gamma_{nj}-\gamma_{n+1,r}-\frac{i\hbar}{2})
\,e^{-i\hbar\d_{\gamma_{nj}}}, \\
\cc_{n}(\la)=-\sum_{j=1}^{n}\prod_{s\neq j}
\frac{\la-\gamma_{ns}}{\gamma_{nj}-\gamma_{ns}} \,
\prod_{r=1}^{n-1}(\gamma_{nj}-\gamma_{n-1,r}+\frac{i\hbar}{2})
\,e^{i\hbar\d_{\gamma_{nj}}} \,.\ee

Note that we also obtain the representation discussed in Section 2
through the following integral formulas which express the generators
$E_{ij}$ of the Lie algebra $\frak{gl}(N)$ in terms of the Yangian
generators
\be\label{rtt3} E_{n,n+1}=\frac{1}{2\pi \hbar}\oint
e_{n}(\la) d\la\;\;,(n=1,\ldots,N-1)\,,\\
E_{n+1,n}=\frac{1}{2\pi\hbar}\oint f_n (\la)
d\la\;\;,(n=1,\ldots,N-1)\,,\\ E_{nn}=\frac{1}{2\pi\hbar}
\oint\limits k_n
(\la)\frac{d\la}{\la}-\frac{1}{2}(n-1)\;\;,(n=1,\ldots,N)\,. \ee
The non-simple root generators  are defined  recursively as
$E_{jk}=[E_{jm},E_{mk}]$ for $j< m<k$ and $j> m>k$. Here, the
integrands are understood as Laurent series and the contours of
integrations are taken around $\infty$.

Let us finally remark, that there is a direct generalization of
the Yangian to the case of quantum group \cite{NT} and it is
possible extend the results of this section to the quantum group
case.
\subsection{The Toda chains and $R$-matrix formalism}
The quantum Toda chain is one of the popular examples of
integrable system. It is described by the Hamiltonian
\be\label{tch1}
H=\sum_{n=1}^N\Big(\frac{p_n^2}{2}+e^{x_n-x_{n+1}}\Big),
\ee
where $[x_n,p_m]=i\hbar\delta_{nm}$. There are two different ways to fix
the boundary conditions: The choice of $x_{\N+1}=\infty$  corresponds to
the open ($GL(N)$) Toda chain; while  the choice of $x_{\N+1}=x_{1}$
corresponds to the periodic (affine) Toda chain.

The  open and periodic Toda chains can be described, uniformly,
by using the $R-$ matrix formalism
\cite{Skl1}. Introduce the Lax operators
\be\label{tch2}
\hspace{1cm}
L_n(\la)\,=\, \left(\begin{array}{cc}\la-p_n & e^{-x_n}\\ -e^{x_n}
& 0
\end{array}\right),\ \ \ \ n=1,\ldots,N\,,
\ee
satisfying the following commutation relations
\be
\hspace{-0.5cm}
R(\la-\mu)(L_n(\la))\otimes I)(I\otimes L_m(\mu))= (I\otimes
L_m(\mu))(L_n(\la)\otimes I)R(\la-\mu)\delta_{nm}\,,
\ee
with the rational $4\times 4$ $R$-matrix
\be
R(\la)=I\otimes I+\frac{i\hbar}{\la}\,P\,.
\ee
The monodromy matrix
\be\label{tch3}
T_{_N}(\la)=L_{\N}(\la)\ldots L_1(\la):=
\left(\begin{array}{cc}A_{\N}(\la) & B_{\N}(\la)\\ C_{\N}(\la) &
D_{\N}(\la)\end{array}\right)
\ee
satisfies the equation
\be\label{tch4}
R(\la-\mu)(T(\la)\otimes I)(I\otimes T(\mu))=
(I\otimes T(\mu))(T(\la)\otimes I)R(\la-\mu)\,.
\ee
In particular, the following commutation relations hold:
\be
[A_{\N}(\la),A_{\N}(\mu)]=[C_{\N}(\la),C_{\N}(\mu)]=0\,,\\
(\la-\mu+i\hbar)A_{\N}(\mu)C_{\N}(\la)\,=\,
(\la-\mu)C_{\N}(\la)A_{\N}(\mu)+i\hbar A_{\N}(\la)C_{\N}(\mu)\,.
\ee
From (\ref{tch4}) it can be easily shown that the trace of the
monodromy matrix
\be\label{tch8}
\wh t_{\N}(\la)=A_{\N}(\la)+D_{\N}(\la)
\ee
satisfies  $[\wh t(\la),\wh t(\mu)]=0$ and is a generating function for
the Hamiltonians of the periodic Toda chain:
\be\label{tch9}
\wh t_{\N}(\la)=\sum_{k=0}^N(-1)^k\la^{N-k}H_k\,.
\ee
We formulate the spectral problems for periodic Toda chain as follows:
\be\label{tch10}
\wh t_{\N}(\la)\Psi_{\raise-3pt\hbox{$\scriptstyle\!\!\bfit E $}}=
t_{\N}(\la;\bfit E)\Psi_{\raise-3pt\hbox{$\scriptstyle\!\!\bfit E
$}}\,, \ee where \be\label{tch11} t_{\N}(\la;\bfit E)=
\sum\limits_{k=0}^N(-1)^k\la^{N-k}E_k\,.
\ee
On the other hand, the operator
\be
A_\N(\la):=\sum_{k=0}^N(-1)^k\la^{N-k}h_k
\ee
is a generating function of the Hamiltonians $h_k$ of the
$N$ particles open Toda chain. The generating functions for the
open Toda chains are connected by the  recursive relations:
\be
A_\N(\la)=(\la-p_\N)A_{\N-1}(\la)+e^{-x_N }C_{\N-1}(\la)\,,\\
C_\N(\la)=-e^{x_N}A_{\N-1}(\la)\,.
\ee

\subsection{The spectral problem for the open Toda chain}
The main goal of this subsection is to apply the results from the
section 2,  to  solution of the spectral problem of the open Toda
chain.

Denote ${\bfit x} =(x_1,\ldots,x_\N)$. To solve the spectral
problem we should find the common eigenfunction of the system of
differential-difference equations:
\be\label{tch16}
A_{{\N}}(\lambda)\psi_{{\bgamma}_{N}}(\bfit x)=
\prod_{m=1}^{N}(\la-\gamma_{\N m})\, \psi_{\bgamma_{N}}(\bfit x)\,,\\
A_{\N-1}(\gamma_{\N j})\psi_{\bgamma_\N}({\bfit x})=i^{1-N}e^{-x_\N}
e^{-i\hbar\d_{\gamma_{Nj}}}\psi_{\bgamma_\N}({\bfit
x})\,,
\ee
$j=1,\ldots,N$. It is worth mentioning that the system (\ref{tch16})
is the quantum counterpart of the Flashka and
McLaughlin \cite{FMc} Darboux transform to separated
variables $(p,x )\rightarrow(\gamma,\theta)$. For the first time the system
(\ref{tch16}) was introduced and solved in the framework of QISM (\cite{KL}).
Below we describe the representation theory solution of the system
(\ref{tch16}).

Let $W'$ and $W$ be the dual irreducible Whittaker modules and
$w'_\N\in W'$, $w_\N\in W$ be the corresponding cyclic
Whittaker vectors. The representation of the Cartan subalgebra is
integrated to the action of the Cartan torus, so  the following
function is well defined \be\label{pair}
\psi_{\gamma_{\N1},\ldots,\gamma_{\N\N}}= e^{-\bfit
x\cdot\brho^{(N)}} \langle
w'_\N,e^{-\sum\limits_{k=1}^Nx_kE_{kk}}w_\N \rangle\,,
\ee
where ${\bfit x}\cdot\brho^{(\N)}$ is the standard product in $\RR^N$.
\begin{de}
The radial projections $A_n (\lambda)$ of the generation functions
 ${\cal A}_n (\lambda)$
(\ref{cas1}) of the central elements of ${\cal U}(\frak{gl}(N))$
are defined by \be\label{h2}
A_n(\la)\psi_{\gamma_{\N1},\ldots,\gamma_{\N\N}} =e^{-{\bfit
x}\cdot\brho^{(N)}} \langle
{w}'_{\N},e^{-\sum\limits_{k=1}^Nx_kE_{kk}} {\cal
A}_n(\la-{\textstyle\frac{i(N-n)\hbar}{2}}) w_{\N} \rangle. \ee
\end{de}
There is the relation between the operators $A_n (\lambda)$ of
different levels: \be\label{rec3}
A_n(\la)=(\la\!-\!p_{n})A_{n-1}(\la)-
e^{x_{n-1}-x_n}A_{n-2}(\la)\,, \ee where $n=1,\ldots,N$ and
$A_{-1}=0,\,A_0=1$. Therefore the $A_n (\lambda)$ is the
generation function of the Hamiltonians of the $n$\ - particles open
Toda chain.

The following theorem identifies  our construction of the matrix
element (\ref{pair}) with the integral formula for the
eigenfunction of the open Toda chain in terms of the Mellin-Barnes
integrals \cite{KL}.
\begin{te}
The matrix element (\ref{pair}) satisfies the set of equations (\ref{tch16}).
\end{te}

It remains to express the matrix element (\ref{pair}) in the
integral form. Substituting the expressions (\ref{wv}), (\ref{wv'}),
and (\ref{m1a}) into (\ref{pair}), we  obtain
\be\label{wf1}
\psi_{\bgamma_{\N}}({\bfit x})\;=e^{-\bfit x\cdot\brho^{(N)}}\\
\hspace{-1.5cm}\times\!\!
\int\limits_{\RR^{\frac{N(N-1)}{2}}} \prod_{n=1}^{N-1}
\frac{\prod\limits_{k=1}^n\prod\limits_{m=1}^{n+1}
\hbar^{\frac{\gamma_{nk}-\gamma_{n+1,m}}{i\hbar}+\frac{1}{2}}\;
\Gamma(\frac{\gamma_{nk}-\gamma_{n+1,m}}{i\hbar}+\frac{1}{2})}
{\prod\limits_{s<p}\left|
\Gamma(\frac{\gamma_{ns}-\gamma_{np}}{i\hbar})\right|^2}
\;e^{\frac{i}{\hbar}
\sum\limits_{n,j=1}^N(\gamma_{nj}-\gamma_{n-1,j})x_n}
\prod_{\stackrel{\scriptstyle n=1}{j\leq n}}^{N-1}d\gamma_{nj}
\,,\hspace{-1.5cm}
\ee
where by definition $\gamma_{nj}=0$ for $j>n$.

In the study of the analytic properties of this solution with
respect to $\bgamma_\N$, it is useful to transform (\ref{wf1}).
Let us change the variables of integration in (\ref{wf1}):
\be\label{wf2}
\hspace{3cm}
\gamma_{nj}\;\to\;\gamma_{nj}-\frac{i\hbar}{n}
\sum_{s=1}^n\rho^{(\N)}_s,\hspace{1cm}(n=1,\ldots,N-1)\,.
\ee
After the change of variables (\ref{wf2}) we shift the domain of
integration $\RR^{\frac{N(N-1)}{2}}$ to the complex plane in such
a way that the domain of integration over the variables
$\gamma_{n-1,j}$ lies above the domain of integration over the
variables $\gamma_{nj}$. Thus, we arrive at the analytic
continuation equal to:
\be\label{wf6}
\psi_{\bgamma_{\N}}(x_1,\ldots,x_\N)\\
=\;\int\limits_{\cal C}\prod_{n=1}^{N-1}
\frac{\prod\limits_{k=1}^n\prod\limits_{m=1}^{n+1}
\hbar^{\frac{\gamma_{nk}-\gamma_{n+1,m}}{i\hbar}}\;
\Gamma(\frac{\gamma_{nk}-\gamma_{n+1,m}}{i\hbar})}
{\prod\limits_{s\neq p}
\Gamma(\frac{\gamma_{ns}-\gamma_{np}}{i\hbar})}
\;e^{\frac{i}{\hbar}
\sum\limits_{n,j=1}^N(\gamma_{nj}-\gamma_{n-1,j})x_n}
\prod_{\stackrel{\scriptstyle n=1}{j\leq n}}^{N-1}d\gamma_{nj}\,,
\ee
where the domain of integration ${\cal C}$ is defined by
the conditions $$\!\min_{j}\{{\rm Im}\,\gamma_{kj}\}>
\max_m\{{\rm Im}\,\gamma_{k+1,m}\}\!$$ for all $ k=1,\ldots,N-1$.
The integral (\ref{wf6}) converges absolutely. Thus, we obtain the integral
representation \cite{KL} for the eigenfunction of  open Toda chain
by purely representation theory methods.

Finally, let us note that the Gelfand-Zetlin type representation
may be generalized to the case of the $Y({\frak g})$, with
${\frak g}$ being an arbitrary simple Lie algebra. This provides the
uniform approach to the solution of the various integrable systems
based on various (quantum) Lie groups. We are planning to discuss
these results elsewhere.

\end{document}